\input amstex
\documentstyle{amsppt}
\magnification=\magstep1

\pageheight{9.0truein}
\pagewidth{6.5truein}

\NoBlackBoxes
\TagsAsMath
\TagsOnLeft

\input xy
\xyoption{matrix}\xyoption{arrow}\xyoption{curve}

\def\edge{\ar@{-}}

\long\def\ignore#1{#1}

\def\la{\Lambda}
\def\lamod{\Lambda{}\operatorname{-mod}}
\def\modla{\operatorname{mod-}{}\Lambda}

\def\Rmod{R{}\operatorname{-Mod}}
\def\Cmod{C{}\operatorname{-mod}}
\def\Rind{R{}\operatorname{-ind}}

\def\Zmod{{\Bbb Z}{}\operatorname{-Mod}}

\def\im{\operatorname{Im}}
\def\id{\operatorname{id}}
\def\length{\operatorname{length}}
\def\soc{\operatorname{soc}}
\def\endosoc{\operatorname{endosoc}}

\def\Hom{\operatorname{Hom}}
\def\End{\operatorname{End}}

\def\add{\operatorname{add}}

\def\inj{\operatorname{in}}
\def\pr{\operatorname{pr}}

\def\NN{{\Bbb N}}
\def\ZZ{{\Bbb Z}}

\def\A{{\Cal A}}

\def\AngHue{{\bf 1}}
\def\Aus{{\bf 2}}
\def\AuSm{{\bf 3}}
\def\Bjork{{\bf 4}}
\def\Coz{{\bf 5}}
\def\CrawBoeI{{\bf 6}}
\def\CrawBoeII{{\bf 7}}
\def\Dung{{\bf 8}}
\def\GruJen{{\bf 9}}
\def\GruJenfinal{{\bf 10}}
\def\Har{{\bf 11}}
\def\HarSai{{\bf 12}}
\def\Huistrong{{\bf 13}}
\def\Huibiel{{\bf 14}}
\def\HuiZimI{{\bf 15}}
\def\HuiZimexchange{{\bf 16}}
\def\HuiZimII{{\bf 17}}
\def\KraSao{{\bf 18}}
\def\Len{{\bf 19}}
\def\Prest{{\bf 20}}
\def\Ring{{\bf 21}}
\def\Ste{{\bf 22}}
\def\Zimhab{{\bf 23}}
\def\Zim{{\bf 24}}
\def\Zimnote{{\bf 25}}

\topmatter

\title  Direct sums of representations as modules over their endomorphism rings 
\endtitle

\rightheadtext{ Endo-structure of direct sums }

\author Birge Huisgen-Zimmermann and Manuel Saor\'\i n \endauthor

\address Department of Mathematics, University of California, Santa Barbara, CA 93106,
USA\endaddress

\email birge\@math.ucsb.edu\endemail

\address Departamento de M\'atematicas, Universidad de Murcia, 30100 Espinardo-MU, Spain
\endaddress

\email msaorinc\@fcu.um.es \endemail

\thanks While carrying out this project, the first-named author was partially supported by an
NSF grant, and the second-named author by grants from the DGES of Spain and the Fundaci\'on
`S\'eneca' of Murcia.
\endthanks
\endtopmatter

\document

\head  1. Introduction, terminology, and background \endhead

Our investigation into the endo-structure of infinite direct sums
$\bigoplus_{i \in I} M_i$ of indecomposable modules $M_i$  --  over a ring $R$ with identity 
--  is centered on the following question:  If $S = \End_R \bigl( \bigoplus_{i \in I} M_i
\bigr)$, how much pressure, in terms of the $S$-structure of
$\bigoplus_{i \in I} M_i$, is required to force the
$M_i$ into finitely many isomorphism classes?  One of the consequences of our principal result
in this direction (Theorem H of Section 4) is as follows.  If all of the $M_i$
are endofinite (think, for instance, of finitely generated or generic modules over an Artin
algebra) and if $(M_t)_{t \in T}$ is a transversal of the isomorphism types of the $M_i$, then
the following conditions (1)-(4) are equivalent:  (1) $T$ is finite;  (2) $\bigoplus_{i
\in I} M_i$ is endo-artinian and $(M_t)_{t \in T}$ is left semi-T-nilpotent; (3) $\bigoplus_{i
\in I} M_i$ is endo-noetherian and $(M_t)_{t \in T}$ is right semi-T-nilpotent; (4) $(M_t)_{t
\in T}$ is left and right semi-T-nilpotent and, for any cofinite subset $T' \subseteq T$, the
endosocle of
$\bigoplus_{t \in T'} M_t$ has finite support in $T'$.  Here we call a family
$(M_i)_{i \in I}$ {\it left semi-T-nilpotent} in case, for every sequence $(i_n)$ of distinct
indices in $I$, every sequence of non-isomorphisms $f_n \in \Hom_R(M_{i_{n+1}}, M_{i_n})$ and
every finitely cogenerated factor module $M_{i_1}/X$ of $M_{i_1}$, the image of $f_1 \cdots
f_n$ is contained in $X$ for $n \gg 0$.  The familiar dual condition of right
semi-T-nilpotence, going back to Harada (see \cite{\Har}), links our results to well-known
theorems addressing the exchange property of direct sums.    

The mentioned equivalence applies with particular strength to Artin algebras, one of the
reasons for this being the following asset (Proposition L): Given any family $(M_i)$ of
finitely generated representations of an Artin algebra, the direct sum
$\bigoplus_{i \in I} M_i$ is endo-artinian if and only if it is $\Sigma$-algebraically
compact (= $\Sigma$-pure injective).  

Our study of $\Sigma$-algebraically compact and, more restrictively,
endo-artinian direct sums hinges crucially on an analysis of endosocles.  The importance of
this invariant in measuring the supply of maps among the $M_i$ emerges clearly in the following
consequence of Theorem H:  A finite dimensional algebra over an algebraically closed field has
finite representation type if and only if the endosocles of all direct sums of indecomposable
left representations of constant finite dimension have finite supports (Corollary N). 

Another point of independent interest lies in the general connections between
$T$-nil\-po\-tence conditions and the endo-structure of direct sums exhibited in Proposition E
of Section 3:  A family $(M_i)_{i \in I}$ of indecomposable modules is right $T$-nilpotent if
and only if
$\bigoplus_{i \in I} M_i$ has the descending chain condition for finitely generated
endo-submodules;  moreover, these conditions are equivalent to the requirement that the direct
sum $\bigoplus_{i \in I} M_i$ be semi-artinian over its endomorphism ring, i.e., that all
endo-factor modules have essential socles.  As a consequence, each $\Sigma$-algebraically
compact module $M$ is a direct sum of indecomposables $M_i$, $i \in I$, with the property that
the family $(M_i)_{i \in I}$ is right
$T$-nilpotent.    

We add a few comments on the background of our project for motivation. It is well known that
the structure of non-finitely generated representations of
$R$, viewed as modules over their endomorphism rings, is decisive in understanding the
behavior of the finitely generated representations.  We point to a few specific
instances of such connections.  As was first observed in
\cite{\Prest} and \cite{\HuiZimII}, finite representation type of $R$ occurs precisely when
all (left) $R$-modules have finite lengths over their endomorphism rings.  Moreover,
the rings with vanishing left pure global dimension are characterized by certain 
endo-chain conditions satisfied by their modules.  Along a different line,
Crawley-Boevey proved that, for tame finite dimensional algebras over an algebraically closed
field, the {\it generic modules\/} (that is, the infinite dimensional indecomposable endofinite
modules) govern the one-parameter families
of finitely generated indecomposable representations in an extremely strong sense
\cite{\CrawBoeI, \CrawBoeII}.  (It appears plausible that the endosocles of direct sums of
generic modules should reflect domestic versus non-domestic representation type; see, e.g.,
Example C(3), due to Ringel
\cite{\Ring}.)  There is a common skein tying the listed scenarios at least loosely together,
namely the fact that
$\Sigma$-algebraic compactness of an
$R$-module is tantamount to the descending chain condition for a selection of its
endo-submodules.  Other results in a related vein link endo-chain conditions of infinite
direct sums of finitely generated representations of Artin algebras to the existence of
preprojective or preinjective partitions, the existence of almost split maps, and the
structure of direct products (see, e.g.,
\cite{\AuSm,
\Huistrong,
\Dung, \KraSao, \AngHue} and our concluding remarks).  These multiple bridges between
endo-chain conditions on one hand, and purely representation-theoretic assets of classes of
modules on the other, motivate our present investigation.    

Section 2 is devoted to exploring {\it endosocle series} of direct sums.  In Section 3, we
compare T-nilpotence conditions for families of modules with endo-chain conditions for their
direct sums, and in Section 4, we follow with our main results, applications, and examples. 
\medskip

\noindent{\it Acknowledgement}  

The authors would like to thank the referee for his numerous helpful suggestions.
\medskip

\noindent{\it Prerequisites}

Recall from \cite{\Zimhab} and \cite{\Zim} that a $p$-{\it functor\/} on $\Rmod$ is a
subfunctor of the forgetful functor $\Rmod \rightarrow \Zmod$ which commutes with direct
products; here $\Rmod$ stands for the category of all left $R$-modules.  Special instances of
$p$-functors can be described as follows:  A {\it pointed matrix\/} over $R$ consists of a
row-finite matrix
$\A = (a_{ij})_{i \in I, j \in J}$ of elements in $R$, paired with a column index
$\alpha \in J$.  Given a pointed matrix $(\A,\alpha)$, we call the following
$p$-functor $[\A,\alpha]$ on $\Rmod$ a {\it matrix-functor\/}:  For any
$R$-module $M$, the subgroup $[\A, \alpha]M$ is defined to be the $\alpha$-th projection of
the solution set in $M^J$ of the homogeneous system
$$\sum_{j \in J} a_{ij} X_j = 0 \qquad \qquad \text{for all}\ \ i \in I.$$ In other words, 
$$[\A, \alpha]M = \{m \in M \mid \exists\ \text{a solution}\ (m_j) \in M^J
\text{\ of the above system with}\ m_{\alpha} = m\}.$$ 
Following \cite{\Zim}, we call $[\A,\alpha]M$ a {\it  matrix
subgroup\/} of $M$, a {\it finite matrix subgroup\/} in case $I$ and $J$ are finite.  We refer
the reader to \cite{\Zim} or \cite{\Huibiel} for more information on
$p$-functors. 

One of the most salient reasons for our present interest in $p$-functors lies in the following
equivalent description of  $\Sigma$-al\-ge\-bra\-i\-cal\-ly compact modules, i.e., of the
modules
$M$ with the property that all direct sums of copies of $M$ are algebraically compact (= pure
injective):  Namely, $M$ is $\Sigma$-al\-ge\-bra\-i\-cal\-ly compact if and only if every
descending chain of p-functorial subgroups of $M$ becomes stationary, a condition which is in
turn equivalent to the DCC for finite matrix subgroups of
$M$ (see \cite{\Zimhab, \GruJen, \Zim, \GruJenfinal});  here we call a subgroup $U$ of the
abelian group underlying $M$ {\it $p$-functorial\/} in case there exists a $p$-functor
$P: \Rmod \rightarrow \Zmod$ with $P(M) = U$.  Since, clearly, every $p$-functorial subgroup
of $M$ is an endo-submodule, this is a well-known instance of a link between the $R$- and
endo-structures of
$M$.   While it guarantees that  endo-artinian modules are
$\Sigma$-al\-ge\-bra\-i\-cal\-ly compact, the converse fails in general; see, e.g.,
\cite{\Coz} or
\cite{\HuiZimI, Theorem 5}.  By contrast,
$M$ is endo-noetherian if and only if $M$ has the ascending chain condition for
$p$-functorial subgroups, since every finitely generated endo-submodule of $M$ is a matrix
subgroup.  In case $M$ is a direct sum of finitely presented modules, the ACC for {\it
finite} matrix subgroups  --  not equivalent to the ACC for arbitrary matrix subgroups in
general  --  already suffices to guarantee endo-noetherianness (Observation 8 of
\cite{\HuiZimII}).   In light of the fact that $p$-functors automatically commute with direct
sums (as arbitrary subfunctors of the forgetful functor do), one deduces
the following observation which will prove useful in connection with dualities.

\example{\it Observation}  For any family $(M_i)_{i \in I}$ of finitely presented modules,
the following conditions are equivalent:

(1) $\bigoplus_{i \in I} M_i$ is endo-noetherian.

(2) $\bigoplus_{i \in I} M_i$ satisfies the ACC for (finite) matrix subgorups.

(3) $\prod_{i \in I} M_i$ satisfies the ACC for (finite) matrix subgroups.

(4) $\prod_{i \in I} M_i$ is endo-noetherian.
\endexample

\head 2. Endosocle series of direct sums of indecomposable modules \endhead

The following elementary observations will be very useful to us in the sequel. 

\proclaim{Lemma A}  Every $\Sigma$-algebraically compact $R$-module $M$ is semi-artinian
over its endomorphism ring $S$, meaning that every $S$-factor module of $M$ has essential
socle.
\endproclaim

\demo{Proof}
 Suppose that $M$ is a $\Sigma$-algebraically compact $R$-module.  As we pointed out above,
$M$ then satisfies the descending chain condition for matrix subgroups which, in particular,
yields the DCC for finitely generated $S$-submodules of $M$.  But this latter condition
clearly entails our claim:  Indeed, given $S$-submodules $U \subsetneqq V$ of $M$, there
exists a minimal
$S$-submodule $W$ of $V$ not contained in $U$, and  any such choice of $W$ gives rise to a
simple submodule $(W+U)/U \subseteq V/U$. 
 \qed 
\enddemo

The converse of Observation A fails in general:  Let R be any left perfect ring which fails to
be left $\Sigma$-algebraically compact (for examples see
\cite{\Zimnote});  then all left $R^{\text{op}}$-modules have essential socles, and therefore
the regular left
$R$-module is semi-artinian over its endomorphism ring.

Recall that the ascending socle series of an
$S$-module $M$ is defined as follows: $\soc_0 M = 0$, $\soc_1 M = \soc_S M$, and, if
$\alpha =
\beta +1$ is a successor ordinal, $\soc_{\alpha} M$ is the preimage of $\soc_S(M /
\soc_{\beta}M)$ under the canonical epimorphism $M \rightarrow M /
\soc_{\beta} M$; for a limit ordinal $\alpha$, finally, $\soc_{\alpha} M$ is defined to be the
union $\bigcup_{\beta < \alpha} \soc_{\beta} M$.  Accordingly, the socle length of $M$ is the
least ordinal $\mu$ with the property that
$\soc_{\mu +1} M = \soc_{\mu} M$.  

Given any left $R$-module $M$ with endomorphism ring
$S$, we will  refer to the $S$-socle of
$M$ as the {\it endosocle\/}, denoted  by $\endosoc M$, and to the
$S$-socle series as the {\it endosocle series\/}; a generic   term of this series will be
labeled 
$\endosoc_{\alpha} M$. The {\it endosocle length} of $M$, finally, is the length of its
$S$-socle series. Whenever we speak of $\endosoc M$ as having a property (X), we are referring
to the
$S$-structure, not the $R$-structure.

Since indecomposable algebraically compact modules have local endomorphism rings
(\cite{\HuiZimI}), Lemma A yields 

\proclaim{Corollary A$'$}  Suppose that $M$ is a $\Sigma$-algebraically compact
$R$-module with endosocle length
$\mu$.  Then
$\endosoc_{\mu} M = M$. 

Hence: If the module $M$ is moreover indecomposable with endomorphism ring $S$, it contains an
exhaustive, well-ordered chain of $S$-submodules with consecutive factors of the form
$(S/J(S))^{(A)}$ (that is, a chain $(M_{\alpha})_{\alpha < \tau}$ of $S$-submodules, where
$\tau$ is an ordinal, with $M_{\alpha} \subseteq M_{\beta}$ for $\alpha < \beta <
\tau$ such that $\bigcup_{\alpha < \tau} M_{\alpha} = M$ and all consecutive factors $M_{\alpha
+ 1} / M_{\alpha}$ are direct sums of copies of $S / J(S)$).  
\quad$\square$
\endproclaim 

In contrast to the preceding observation, {\it direct sums\/} of $\Sigma$-algebraically
compact modules may have trivial endosocles  --  think of $M =
\bigoplus_{n \in \NN} \ZZ/ (p^n)$, for instance, where $p$ is  prime.  This can be seen as a
pronounced failure of the direct sum of $\Sigma$-algebraically compact modules to inherit
 crucial properties implied by
$\Sigma$-algebraic compactness.  

Recall that any $\Sigma$-algebraically compact module $M$ is a direct sum of submodules with
local endomorphism rings.  Our main interest being in families $(M_i)$ of
$\Sigma$-algebraically compact or, more restrictively, endo-artinian modules, most of our
discussion will therefore focus on the situation where all of the
$M_i$ have local endomorphism rings.   In this setting, we record how the endosocle of
$\bigoplus_{i \in I} M_i$ relates to the endosocles of the individual $M_i$.    In case all of
the $M_i$ are pairwise isomorphic, the connection is straightforward, otherwise it reflects
the behavior of maps among the summands $M_i$ as follows.

\proclaim{Lemma B}  Let $(M_i)_{i \in I}$ be a family of modules with local
endomorphism rings, and $S = \End_R(\bigoplus_{i \in I} M_i)$.
\medskip

{\rm (1)}  First suppose that $M_i = N$ for all $i \in I$.  Then:
\roster
\item"$\bullet$"  
$\endosoc \bigl( N^{(I)}\bigr) = \bigl(\endosoc(N) \bigr)^{(I)}$, and the endosocle of
$N^{(I)}$ consists of a single homogeneous component.
\item"$\bullet$"  The endosocle of $N^{(I)}$ is finitely generated as an endo-submodule if and
only if the same is true for the endosocle of $N$.
\endroster

{\rm (2)}  In general, $\endosoc \bigl(\bigoplus_{i \in I} M_i\bigr) = \bigoplus_{i \in I}
B_i$, where each $B_i$ is the $\End_R(M_i)$-submodule of $M_i$ consisting  precisely of those
elements of $M_i$ which are annihilated by all non-isomorphisms in
$\bigcup_{j \in I} \Hom_R(M_i, M_j)$.

In particular:  If $M_i \not\cong M_j$ for $i \ne j$, then each $B_i$ is a semisimple
$\End_R(M_i)$-submodule of $M_i$, as well as a semisimple $S$-submodule of the direct sum
$\bigoplus_{i
\in I} M_i$.
\smallskip

{\rm (3)} Whenever $(M_n)_{n \in \NN}$  is a family of pairwise nonisomorphic modules with
local endomorphism rings which permits a sequence of monomorphisms 
$M_{1}
\hookrightarrow  M_{2} \hookrightarrow M_{3} \hookrightarrow \cdots$, the direct sum
$\bigoplus_{n
\in
\NN} M_{n}$ has trivial endosocle.
  
\endproclaim

\demo{Proof} Set $M = \bigoplus_{i \in I} M_i$, and observe that each $S$-submodule $U$ of $M$
is of  the form $U = \bigoplus_{i \in I} (U \cap M_i)$, where the intersection $U \cap M_i$ is
an $\End_R(M_i)$-submodule of $M_i$.  We deduce  that each simple $S$-submodule of $M$ can be
written in the form $S x_i$ for a nonzero element $x_i$ of some
$M_i$.

Part (1) is immediate from these remarks. 

(2)  Set $B_i = (\endosoc M) \cap M_i$.  To see that $B_i$ is contained in the annihilator of
the set of non-isomorphisms in $\bigcup_{j \in I} \Hom_R(M_i, M_j)$, observe that, due to the
locality of the rings
$\End_R(M_i)$, each of the indicated non-isomorphisms belongs to the Jacobson radical of $S$.

For the other inclusion, fix $i
\in I$, set $I(i) = \{j \in I \mid M_j \cong M_i \}$, and $S(i) = \End_R \bigl( \bigoplus_{j
\in I(i)} M_j \bigr)$.  If $x \in M_i$ is annihilated by all non-isomorphisms in the above
union, then
$x$ clearly belongs to the $\End_R(M_i)$-socle of $M_i$, and consequently $S(i) x$ is
contained in the $S(i)$-socle of $\bigoplus_{j \in I(i)} M_j$ by part (1).  Since $Sx =
S(i)x$, this shows that
$Sx \subseteq \endosoc M$, which yields $x \in B_i$ as required. 

The final assertion under (2) is now obvious, as is (3).         
\qed \enddemo

\definition{Examples C}

{\rm (1) Let $R = \ZZ$ and $p$ a prime.  The terms of the endosocle series of $M = \bigoplus_{n
\in \NN} \ZZ/(p^n)
\oplus \ZZ(p^{\infty})$ are as follows:   $\endosoc_{\alpha} M$ equals the copy of $\ZZ/
(p^{\alpha})$ inside $\ZZ(p^{\infty})$ for $\alpha < \omega$, and $\endosoc_{\alpha} M =
\ZZ(p^{\infty})$ for
$\alpha \ge \omega$.   
\smallskip 

{\rm (2)}  Let $\la$ be the Kronecker algebra. If $(M_n)_{n \in \NN}$ is the family of all
preprojective modules in
$\Rind$, then
$\endosoc \bigl( \bigoplus_{n \in \NN} M_n \bigr) = 0$; in fact, $\endosoc \bigl( \bigoplus_{n
\ge m} M_n \bigr) = 0$ for all $m$.   On the other hand, if $M' =
\bigoplus_{n \in \NN} M'_n$ is the direct sum of all preinjective modules, with
$\length(M'_n) = 2n - 1$,  the endosocle of
$M'$ equals $M'_1 \oplus \soc_{\la} M_2'$, and, for $m \ge 2$, we have $\endosoc
\bigl(\bigoplus_{n
\ge m} M'_n\bigr) = M'_m$.

{\rm (3)}  The following is due to Ringel \cite{\Ring}:  If $\la$ is a string algebra (i.e., a
finite dimensional monomial relation algebra which is special biserial), but fails to be of
domestic representation type, then there exists a family
$(M_n)_{n \in \NN}$ of pairwise non-isomorphic generic $\la$-modules allowing for consecutive
embeddings
$M_1 \hookrightarrow M_2 \hookrightarrow \dots$; by Lemma B(3), we infer that the
endosocle of the direct sum of the $M_n$'s is zero.   
\enddefinition        

We continue to assume that $M =
\bigoplus_{i \in I} M_i$, where all $M_i$ have local endomorphism rings.  In exploring the
endo-structure of $M$, a different sequence of iterated endosocles  --  not forming an
ascending chain  --  is frequently more helpful than the series we just discussed.  This is
due to the fact that the traditional endosocle series may be infinite and still get stalled 
within a finite subsum of
$M_i$'s  (see Example C(1)).  By contrast, any infinite sequence of nonzero terms of the
following {\it relative endosocle series} involves infinitely many summands $M_i$.  This often
makes it a more effective tool in studying homomorphisms among non-isomorphic $M_i$'s.  

 As usual, we start with $\endosoc'_0  M = 0$;  moreover, we set $I_0 =
\varnothing$.  Next, we set 
$\endosoc'_1 M = \endosoc M$, and denote by $I_1$ the support of $\endosoc'_1 M$ in $I$. 
Assume that the terms $\endosoc'_{\beta} M$ have already been defined for all ordinal numbers
$\beta < \alpha$.  In case
$\alpha$ is a limit ordinal, we set $\endosoc'_{\alpha} M  = 0$.  If, on the other hand,
$\alpha$ is a successor ordinal, say $\alpha = \beta +1$,  we let $I_{\beta}$ be the support of
$\sum_{\gamma \le \beta} \endosoc'_{\gamma} M$, and define  
$\endosoc'_{\alpha} M$ to be the endosocle of the trimmed direct sum $\bigoplus_{I
\setminus I_{\beta}} M_i$.  In particular, we see that $I_{\alpha}$ is defined for each
ordinal $\alpha$ and equals the support of $\sum_{\beta \le \alpha} \endosoc'_{\beta} M =
\bigoplus_{\beta \le \alpha} \endosoc'_{\beta} M$. Finally, we refer to the least ordinal
number
$\mu$ such that $\endosoc'_{\mu + 1} M = 0$ (or, equivalently, the least ordinal number $\mu$
with $I_{\mu+1} = I_{\mu}$)  as the {\it relative endosocle length} of $M$.

Clearly,  the isomorphism types of the $R$-submodules $\endosoc'_{\alpha} M$ are isomorphism
invariants of
$M$, and the sum $\sum_{\alpha \le \mu} \endosoc'_{\alpha} M$ is direct by construction.  The
second of the following observations will be used repeatedly in the next section.  Both are
straightforward from the definitions.

\proclaim{Lemma D} Let $M = \bigoplus_{i \in I} M_i$, where each $M_i$ has local
endomorphism ring.  

{\rm (1)}  If each nontrivial subsum of $M_i$'s has nontrivial endosocle, then $I_{\mu} = I$,
where $\mu$ is the relative endosocle length of $M$.  

{\rm (2)}  Let $\alpha$ be a successor ordinal, say $\alpha = \beta + 1$, such that
$\bigoplus_{i \in I \setminus I_{\beta}} M_i$ has essential endosocle.  If 
$U$ is an endo-submodule of $M$ which is not contained in
$\bigoplus_{i \in I_{\alpha}} M_i$, then $\bigl( J(S)U
\bigr) \cap \endosoc'_{\alpha} M\ne 0$.  \qed    \endproclaim

Observe that all hypotheses of Lemma D are satisfied if $M$ is $\Sigma$-algebraically
compact.

\head 3. T-nilpotence versus endo-chain conditions
\endhead

Harada \cite{\Har} called a family $(M_i)_{i \in I}$ of indecomposable $R$-modules {\it
semi-T-nilpotent} in case, for any sequence $(i_n)_{n \in \NN}$ of distinct indices in $I$,
any family of non-isomorphisms $f_n \in \Hom_R(M_{i_n}, M_{i_{n+1}})$, and any finitely
generated $R$-submodule $X$ of $M_{i_1}$, there exists $n_0 \in \NN$ with 
$$f_{n_0} f_{n_0 -1} \cdots f_1(X) = 0.$$   Obviously, this is a condition which can be tested
`pointwise', i.e., the family $(M_i)$ satisfies it precisely when, for each family of indices
and non-isomorphisms $f_n$ as above and any element $x \in M_{i_1}$, there exists a number
$n_0$ with
$f_{n_0} f_{n_0 -1} \cdots f_1(x) = 0$.  To make room for a dual concept, we will refer to the
described condition as {\it right semi-T-nilpotence}. The following twin concept of {\it left
semi-T-nilpotence} of a family
$(M_i)_{i \in I}$ of indecomposables requires that, for any sequence
$(i_n)_{n \in \NN}$ of distinct indices in $I$, any family of non-isomorphisms
$f_n \in \Hom_R(M_{i_{n+1}}, M_{i_n})$, and any finitely cogenerated factor module $M_{i_1}/X$
of the $R$-module $M_{i_1}$, there exists a natural number
$n_0$ such that 
$$\im(f_1 f_2 \cdots f_{n_0}) \subseteq X.$$   Left semi-T-nilpotence clearly implies the
following elementwise condition:  For each family of indices and homomorphisms
$f_n$ as specified above, and for each nonzero $y \in M_{i_1}$, there exists $n_0 \in \NN$ such
that
$y \notin \im(f_1 f_2 \cdots f_{n_0})$.  While the converse fails in general, it does hold if
all of the $M_i$'s are finitely cogenerated  --  so, in particular, if they have finite length 
--  in which case either condition is equivalent to the requirement that the family of
homomorphisms among the $M_i$'s be artinian in the sense of Auslander (see below). 
 
The stronger conditions of {\it  right/left T-nilpotence}, finally, call for the same
conclusions as the corresponding semi-T-nilpotence conditions, but on waiving the premise that
the families of indices considered be free of repetitions.  

 We point out that, for the special case of finitely generated indecomposable modules $M_i$
over an Artin algebra, Auslander introduced alternate terms for these nilpotence conditions
(see
\cite{\Aus}).  He labeled a family of morphisms among such modules noetherian (resp.,
artinian) in case, for every countable sequence of  non-isomorphisms $(f_n)$ selected within
the family  --  consecutively composable in the appropriate sense  --  there exists an index
$n_0$ with
$f_{n_0} \cdots f_1 = 0$ (resp., $f_1 \cdots f_{n_0} = 0$).  So, in this situation, the family
$(M_i)$ is right semi-T-nilpotent if and only if it is right T-nilpotent, if and only if the
family of homomorphisms among the
$M_i$ is noetherian.  (We alert the reader to the fact that this condition automatically comes
with endo-{\it artinian\/}ness of $\bigoplus_{i \in I} M_i$  --  see Proposition E below  --  
but may fail for endo-noetherian direct sums $\bigoplus_{i \in I} M_i$.)

The second of the following known characterizations of right semi-T-nilpotent families of
modules will be used freely in the sequel, while the third will enter into an application of
our main theorem.  The implications `(3)$\implies$(2)$\iff$(1)' are due to Harada
(\cite{\Har}, and the remaining one was filled in by W. Zimmermann and the firstnamed author
\cite{\HuiZimexchange}.  Recall that a module
$M$ is said to have the {\it exchange property} in case, for every equality of the form $M'
\oplus A = \bigoplus_{l \in L} A_l$ with $M' \cong M$, there exist submodules
$B_l \subseteq A_l$ such that $M' \oplus A = M' \oplus \bigoplus_{l \in L} B_l$.

\proclaim{Known facts}  Suppose that all $M_i$ have local endomorphism rings.  Then the
following conditions are equivalent:

{\rm (1)}  The family  $(M_i)_{i \in I}$ is right semi-T-nilpotent. 

{\rm (2)}  The Jacobson radical of $\End_R
\bigl(\bigoplus_{i \in I} M_i \bigr)$ coincides with the set of those endomorphisms $f$ which
have the property that all compositions $\pr_j \cdot f \cdot \inj_i$ are non-isomorphisms
(here $\pr_j$ and $\inj_i$ denote the canonical projections and injections, respectively).
   
{\rm (3)}  The direct sum $\bigoplus_{i \in I} M_i$ has the exchange property. \endproclaim

We start by relating right semi-T-nilpotence to the endo-semi-artinian condition we encountered
in Lemma A.  These connections will not only feed into the proof of our main result, but
are of interest as supplements.  Part (3) of the following proposition is a variant of
\cite{\HuiZimII, Proposition 4}.  Moreover, we note that the implications
(a)$\implies$(c)$\implies$(b) hold without the blanket hypothesis under (2).

\proclaim{Proposition E}  Suppose that $(M_i)_{i \in I}$ is a family of indecomposable left
$R$-modules, and set $M = \bigoplus_{i \in I} M_i$. 

{\rm (1)}  If $(M_i)$ is a right T-nilpotent family, then all of the $M_i$ have local
endomorphism rings. 

{\rm (2)}  If all $M_i$ have local endomorphism rings, the following conditions are
equivalent: 
\widestnumber\item{(a)\quad}\roster
\item"{\rm (a)}" The family $(M_i)$ is right T-nilpotent.
\item"{\rm (b)}" The direct sum $M$ is semi-artinian over its endomorphism ring.
\item"{\rm (c)}" $M$ satisfies the DCC for cyclic endo-submodules (or, equivalently, the DCC
for finitely generated endo-submodules).  
\endroster

As a consequence, every $\Sigma$-algebraically compact $R$-module $M$ is a direct sum ranging
over a right T-nilpotent family of indecomposable modules.
\smallskip

{\rm (3)} The following condition is sufficient for right semi-T-nilpotence of the family
$(M_i)$:  Each $M_i$ has local endomorphism ring, and every descending chain of matrix
functors on
$\Rmod$ becomes uniformly stationary on almost all of the $M_i$'s. 
\smallskip 

\endproclaim 

\demo{Proof} Set $S= \End_R(M)$.

(1)  Suppose that $M_i$ belongs to a right T-nilpotent family, and let $f \in \End_R(M_i)$ be a
non-isomorphism.  Then $f$ fails to be a monomorphism, for, given any element $x \in M$, there
exists $m \in \NN$ with $f^m(x) = 0$.  To see that, for any $g \in \End_R(M_i)$, the map
$\id_{M_i} - gf$ is an isomorphism, it thus suffices to check injectivity.  To that end, we
observe that $x = gf(x)$ implies $x = (gf)^n(x)$ for all $n \in \NN$, and hence $x=0$ by our
T-nilpotence hypothesis.  

(2) Start by recalling that the descending chain condition for cyclic submodules is equivalent
to that for finitely generated submodules, due to Bj\/ork
\cite{\Bjork}.  The canonical embeddings and projections coming with the direct sum
$\bigoplus_{i
\in I} M_i$ will be denoted by $\inj_i$ and $\pr_i$, respectively.

To back up `(a)
$\implies$ (c)', suppose (a) holds. 

First we show that, for any element $m = (m_i)_{i \in I} \in M$, the cyclic $S$-module $Sm/
J(S)m$ is semisimple.  Indeed, right semi-T-nilpotence of $(M_i)$ guarantees that the Jacobson
radical $J(S)$ of $S$ contains all endomorphisms $f$ with the property that the compositions
$\pr_i f \inj_j$ are non-isomorphisms for arbitrary $i,j \in I$.  We will make the assumption
that $M_i \not\cong M_j$ whenever
$i \ne j$; this does not affect the generality of our argument, for every endo-submodule
$U$ of a direct sum of powers $\bigoplus_{i \in I} M_i^{(L_i)}$ is of the form
$\bigoplus_{i \in I} U_i^{(L_i)}$ for a suitable endo-submodule $U_i$ of
$M_i$. Consequently,
$$Sm/ J(S)m = \bigoplus_{\text{finite}} \End_R(M_i) m_i/ J \bigl(\End_R(M_i) \bigr)m_i,$$
where each of the summands on the right is an $S$-module, the $S$-structure of which coincides
with the pertinent 
$\End_R(M_i)/ J \bigl(\End_R(M_i) \bigr)$-structure.

To the contrary of our claim assume that $M$ contains a strictly descending chain $S x_1
\supsetneqq  S x_2 \supsetneqq S x_3 \supsetneqq \cdots$ of cyclic $S$-submodules.  Now the
induced sequence in the semisimple factor module $S x_1/J(S)x_1$ does become stationary,
yielding an integer $m$ such that $Sx_m \subseteq Sx_{m+1} + J(S)x_1$.  By possibly dropping
some terms from our descending chain, we may assume that $Sx_2 \subseteq Sx_3 + J(S)x_1$.  We
apply the same argument to $Sx_2 / J(S) x_2$ and iterate, whence an obvious induction permits
us to thin out our original sequence of cyclics so as to assure $Sx_n \subseteq Sx_{n+1} +
J(S)x_{n-1}$ for all $n \ge 2$.  For each $n \ge 2$ we thus obtain an equality $x_n =
f_{n+1}(x_{n+1}) + g_{n-1}(x_{n-1})$ with $f_{n+1} \in S$ and $g_{n-1} \in J(S)$.  By
substituting $x_2 = f_3(x_3) + g_1(x_1)$ into the equality $x_3 = f_4(x_4) + g_2(x_2)$, we
deduce that
$(1 - g_2 f_3)(x_3) = f_4(x_4) + g_2 g_1(x_1)$, and multiplication by the inverse of $1 - g_2
f_3$ yields an equality of the form
$x_3 = f_4'(x_4) + h_1(x_1)$ with $h_1 \in J(S)$.  In a next step, we insert this expression
for $x_3$ into $x_4 = f_5(x_5) + g_3(x_3)$, thus obtaining an element $h_2 \in J(S)$ with the
property that $x_4 = f_5'(x_5) + h_2 h_1(x_1)$ for a suitable $f_5' \in S$.  Again we iterate,
and arrive at a sequence of elements $h_1, h_2, h_3, \dots \in J(S)$ satisfying $x_n =
f_{n+1}'(x_{n+1}) + h_{n-2}
\cdots h_2\,h_1(x_1)$ for suitable choices of $f_n' \in S$.  

T-nilpotence of the family $(M_i)_{i\in I}$, combined with Koenig's graph theorem, will now
provide us with a natural number
$N$ such that $h_N \cdots h_2\,h_1(x_1) = 0$.  For, if we had $h_n
\cdots h_2\,h_1(x_1) \ne 0$ for all $n$, the following graph $\Cal G$ would have paths of
arbitrary length:  the roots of $\Cal G$ are the nonzero components $x_{1i} \in M_i$ of the
element $x_1$  --  finite in number, say $x_{1,i(1)}, \dots, x_{1, i(m)}$, and the edge paths
of length $l$ in $\Cal G$ correspond to the nonzero evaluations 
$$\bigl(\pr_{j(l)} h_l \inj_{j(l-1)}\bigr)\bigl(\pr_{j(l-1)} h_{l-1} \inj_{j(l-2)}
\bigr)
\cdots \bigl(\pr_{j(1)} h_1 \inj_{j(0)}\bigr) \bigl(x_{1,j(0)}\bigr),$$  for suitable
indices $j(s) \in I$ such that $j(0) \in
\{i(1), \dots, i(m)\}$.  Koenig's graph theorem would then imply the existence of an infinite
path.  But this is incompatible with our T-nilpotence condition, since each of the maps
$\pr_{j(s)} h_s \inj_{j(s-1)}: M_{j(s-1)} \rightarrow M_{j(s)}$ is a non-isomorphism by
our choice of the $h_n$.    

Thus, we obtain a number $N$ as described and infer
$x_{N+2} \in Sx_{N+3}$, contrary to our assumption that the chain of $Sx_n$ be strictly
descending. 
\medskip

To justify the implication `$(c) \implies (b)$', carry over the argument proving Lemma A.

For `$(b) \implies (a)$', see \cite{\Ste, Chapter VIII}.
\smallskip

(3)  Suppose that $(i_n)_{n \in \NN}$ is a sequence of distinct elements of $I$, and
$f_{i_n} : M_{i_n} \rightarrow M_{i_{n+1}}$ a non-isomorphism for $n \in \NN$.  Given
$x \in M_{i_1}$, consider the following chain of principal $S$-submodules of $M$:
$$S f_{i_1} x \supseteq S f_{i_2} f_{i_1} x \supseteq S f_{i_3} f_{i_2} f_{i_1} x
\supseteq \cdots.$$ All terms of this chain are matrix subgroups of the $R$-module $M$ (see
Section 1), say
$S f_{i_n} \cdots f_{i_1} x = P_n(M)$, where $P_n$ is a matrix functor.  Since the class of
matrix functors on
$\Rmod$ is closed under intersections, we can normalize to the situation where $P_1 \supseteq
P_2 \supseteq P_3 \supseteq \cdots$.  Assuming that this chain becomes uniformly stationary on
amost all of the $M_i$'s, we obtain a natural number $N$, together with a finite subset $I_0
\subseteq I$, such that $P_n \bigl( \bigoplus_{i \in I
\setminus I_0} M_i \bigr) = P_N \bigl( \bigoplus_{i \in I \setminus I_0} M_i\bigr)$ for all
$n \ge N$.  Since the $i_n$ are pairwise different, we can adjust $N$ upwards, if necessary,
so as to guarantee that, moreover, $f_{i_n} \cdots f_{i_1} x
\in \bigoplus_{i \in I \setminus I_0} M_i$ for all $n \ge N$.  Denoting the
$R$-endomorphism ring of $\bigoplus_{i \in I \setminus I_0} M_i$ by $S'$, we infer that the
chain 
$$S' f_{i_n} \cdots f_{i_{N+1}} \bigl( f_{i_N} f_{i_{N-1}} \cdots f_{i_1} x \bigr),
\quad \quad n \ge N$$ is stationary.  We will deduce that $f_{i_N} \cdots f_{i_1} x = 0$.
Indeed, our setup yields a map $g \in S'$ with
$g f_{i_{N+1}} f_{i_N} \cdots f_{i_1} x =  f_{i_N} \cdots f_{i_1} x$.  Clearly, we may assume
that $g$ is a map from $M_{i_{N+2}}$ to $M_{i_{N+1}}$.  If $f_{i_N} \cdots f_{i_1} x$ were
nonzero, Nakayama's Lemma would exclude
$g f_{i_{N+1}}$ from the Jacobson radical of $\End_R(M_{i_N})$, and locality of this
endomorphism ring would guarantee $g f_{i_{N+1}}$ to be a unit of
$\End_R(M_{i_N})$.  But this is impossible, since $M_{i_{N+2}}$ is indecomposable and
$f_{i_{N+1}}$ a non-isomorphism.   Hence $f_{i_N} \cdots f_{i_1} x = 0$ as claimed, and the
argument is complete. \qed \enddemo

Part (2) of the preceding proposition can actually be de-specified from the situation of direct
sums regarded as modules over their endomorphism rings:  If $S$ is any ring and
$M$ a left
$S$-module, then the DCC for finitely generated submodules forces $M$ to be semi-artinian; the
latter condition in turn implies that, given any $m \in M$ and any sequence
$(s_1, s_2, \dots)$ in $J(S)$, we have $s_n \cdots s_1 m = 0$ for $n \gg 0$.  Moreover, these
three conditions are equivalent in case $Sm/J(S)m$ is semisimple for every element $m \in M$.

\head 4. Main results and examples \endhead    

Theorem H below will zero in on our primary concern:  Namely to significantly weaken
conditions on the endo-structure of
$\bigoplus_{i \in I} M_i$ which are known to guarantee that the $M_i$ fall into finitely many
isomorphism classes.  An instance of such a condition is endofiniteness of $\bigoplus_{i \in
I} M_i$; see \cite{\CrawBoeII}.  Our result also provides background for the following
well-known characterization of the Artin algebras
$\la$ having finite representation type, in terms of morphisms among the indecomposable
finitely generated left
$\la$-modules, due to Auslander ({\cite{\Aus}}):  Namely,
$\la$ has finite representation type precisely when every countable family of non-isomorphisms
is both artinian and noetherian in the sense given at the beginning of Section 2.  Related
results for families of finitely generated indecomposable modules were recently obtained by
Dung
\cite{\Dung, Theorem 3.3 and Corollary 3.12}.  We, too, are particularly interested in the
situation where all of the
$M_i$'s have finite length over $R$  --  next to the case where the $M_i$'s are generic  -- 
and will address it in subsequent corollaries and examples.

\proclaim{Lemma F}  Let $(M_i)_{i \in I}$ be a right semi-T-nilpotent family of pairwise
non-isomorphic modules with local endomorphism rings.  If, for each cofinite subset $I'$ of
$I$, the direct sum
$\bigoplus_{i
\in I'}M_i$ is finitely generated over its endomorphism ring, then $I$ is finite.
\endproclaim

\demo{Proof}  Suppose that the endo-condition of our claim is satisfied, and set $M =
\bigoplus_{i \in I} M_i$. For any subset $I' \subseteq I$, we denote by
$\pr_{I'}$ the canonical projection $M \rightarrow \bigoplus_{i \in I'} M_i$ along
$\bigoplus_{i \in I \setminus I'} M_i$, preserving the convention $\pr_i =
\pr_{\{i\}}$ for $i \in I$.  Moreover, we set $S = \End_R(\bigoplus_{i \in I} M_i \bigr)$. 
Keeping in mind that every $S$-submodule $U$ of
$M$ is of the form $\bigoplus_{i \in I} (M_i \cap U)$, where each $M_i \cap U$ is an
$\End_R(M_i)$ submodule of $M_i$, we find that, given any $S$-generating set $G$ of such a
module $U$, the collection $\{\pr_i(x) \mid i \in I, x \in G\}$ of projections again belongs
to $U$ (and thus generates $U$). 

By hypothesis, we can therefore pick a finite family of elements $(x_{l_1})_{l_1
\in L_1}$ in $\bigcup_{i \in I} M_i$ with the property that $M = \sum_{l_1 \in L_1} S
x_{l_1}$;  by $K_1 \subseteq I$ we denote the (finite) union of the supports of the
$x_{l_1}$.  Setting $S_1 = \End_R(\bigoplus_{i \in I \setminus K_1} M_i)$ and applying our
hypothesis to the pared-down direct sum, we see that each of the modules $S_1\pr_{I
\setminus K_1}(S x_{l_1})$ can be written as a finite sum, $S_1 \pr_{I
\setminus K_1}(S x_{l_1}) = \sum_{l_2 \in L_2} S_1 f_{l_2 l_1}(x_{l_1})$, where $L_2$ is a
finite index set disjoint from $L_1$ and the $f_{l_2 l_1}$ are homomorphisms from suitable
summands
$M_i$ of $\bigoplus_{k \in K_1} M_k$ to summands $M_j$ of $\bigoplus_{k \in I\setminus K_1}
M_k$; here the sources and targets of the $f_{l_2 l_1}$ may appear repeatedly.  By allowing
zero maps among the $f_{l_1 l_2}$, we can make the same index set $L_2$ work for all
$l_1 \in L_1$ simultaneously.  We thus obtain 
$$\align M &= \bigoplus_{k \in K_1} M_k \ + \ \sum_{l_1 \in L_1} S_1 \pr_{I\setminus K_1}
(Sx_{l_1})\\
 &= \bigoplus_{k \in K_1} M_k \ + \ \sum_{l_1 \in L_1, l_2 \in L_2}  S_1 f_{l_2 l_1}
(x_{l_1}).\endalign$$

Next we let $K_2 \subseteq I$ be a finite set containing $K_1$ and the supports of the
elements $f_{l_2 l_1}(x_{l_1})$ for all choices of $l_1$, $l_2$.   Now we set $S_2 =
\End_R(\bigoplus_{i \in I \setminus K_2} M_i)$ and, using the same guidelines as for the
choice of
$L_1$ and the
$f_{l_2 l_1}$'s, we choose a finite index set $L_3$, disjoint from $L_1 \cup L_2$, and
homomorphisms
$f_{l_3 l_2}$, each from a suitable summand $M_i$ of $\bigoplus_{k \in K_2} M_k$ to a summand
$M_j$ of $\bigoplus_{k \in I \setminus K_2} M_k$, such that 
$$S_2 \pr_{I \setminus K_2} \bigl( S_1f_{l_2 l_1} (x_{l_1}) \bigr)\ =\ \sum_{l_3 \in L_3} S_2
f_{l_3 l_2} f_{l_2 l_1} (x_{l_1}),$$ for all $l_1 \in L_1$ and $l_2 \in L_2$.  This yields the
equality
$$M \ = \bigoplus_{k \in K_2} M_k \ + \ \sum_{l_1 \in L_1, l_2 \in L_2, l_3 \in L_3} S_2
f_{l_3 l_2} f_{l_2 l_1} (x_{l_1}).$$  

We repeat the above procedure to  inductively obtain pairwise disjoint finite index sets $L_1,
L_2, L_3, \dots$, finite subsets $K_1 \subseteq K_2 \subseteq K_3 \subseteq
\cdots$ of $I$, and homomorphisms $f_{l_{n+1} l_n}$  --  for $l_n \in L_n$ and
$l_{n+1}
\in L_{n+1}$  --  from summands $M_i$ of $\bigoplus_{k \in K_n} M_k$ to summands
$M_j$ of $\bigoplus_{i \in I \setminus K_n} M_i$, respectively, such that  for each $n
\in \NN$,
$$M \ = \bigoplus_{k \in K_n} M_k \ + \ \sum_{l_m \in L_m, 1 \le m \le n} S_n f_{l_{n+1}
l_n}\cdots f_{l_2 l_1} (x_{l_1}),$$   where  $S_n = \End_R(\bigoplus_{i \in I \setminus K_n}
M_i)$.

Let $\Cal G$ be the graph having as vertices the indices in $\bigcup_{i \in \NN} K_i$ and as
edge paths those concatenations $f_{l_{n+1} l_n} \cdots f_{l_2 l_1}$ which have the property
that
$f_{l_{n+1} l_n} \cdots f_{l_2 l_1}(x_{l_1}) \ne 0$.  Since all of the maps $f_{l_{i+1} l_i}$
are non-isomorphisms by construction, our T-nilpotence hypothesis excludes the occurrence of
an infinite path.  Therefore, K\"onig's Graph Theorem guarantees the existence of an upper
bound on the lengths of paths in $\Cal G$   --  N say  --  which shows that $M = \bigoplus_{k
\in K_N} M_k$.  But this says that  $I = K_N$ is finite as required. \qed \enddemo

We now come to our key lemma.  In spite of the dual formats of Lemmas F and G, the latter
cannot be obtained by dualizing the preceding argument.

\proclaim{Lemma G}  Let $(M_i)_{i \in I}$ be a left semi-T-nilpotent family of modules with
local endomorphism rings.  If, for each cofinite subset $I'$ of
$I$, the direct sum $\bigoplus_{i \in I'}M_i$ is finitely cogenerated over its endomorphism
ring, then
$I$ is finite.
\endproclaim

\demo{Proof} Again suppose that the endo-condition of our claim is satisfied, and set $M =
\bigoplus_{i \in I} M_i$, 
$S = \End_R(M)$.  Assume that $I$ is infinite.  In order to deduce that the family $(M_i)_{i
\in I}$ then fails to be left semi-T-nilpotent, we consider the relative endosocle series  
$(\endosoc'_{\alpha} M)_{\alpha}$, as introduced at the end of Section 2.  As before, we
denote the support of $\bigoplus_{\beta \le
\alpha}
\endosoc'_{\beta} M$ by $I_{\alpha}$.  By hypothesis, $\endosoc M$ is finitely generated,
whence the support
$I_1$ is finite by Lemma B(2).  An obvious induction, based on the hypothesis that the
endosocle of any cofinite subsum
$\bigoplus_{i \in I'} M_i$ has finite support in $I'$, thus yields finiteness of
$I_k$ for all $k \in \NN$. 

Our hypothesis moreover guarantees that each of the relative endosocles $\endosoc'_{k+1} M$ is
essential in the corresponding direct sum 
$\bigoplus_{i \in I \setminus I_{k}} M_i$ viewed as a module over its endomorphism ring
$S_k$.  Since, under the canonical embedding of $S_k$ into $S$, we have $J(S_k) \subseteq
J(S)$, we deduce from Lemma D(2) that, for any nonzero $S_k$-submodule $U$ of
$\bigoplus_{i \in I \setminus I_{k}} M_i$, the intersection $\bigl( J(S) U \bigr) \cap
\endosoc'_{k} M$ is nonzero. For each $k \in \NN$, we now consider a descending chain $A_{k1}
\supseteq A_{k2} \supseteq A_{k3}
\supseteq \dots$ of nonzero $S_k$-submodules of $\endosoc'_{k} M$, where the $A_{ki}$ are
inductively defined as follows:  $A_{k1} = (J(S) \endosoc'_{k+1} M) \cap \endosoc'_{k} M$ for
$k 
\ge 1$, and 
$$A_{k,i+1} = (J(S) A_{k+1,i})\cap \endosoc'_k M.$$  Being contained in the finitely generated
semisimple $S_{k-1}$-module $\endosoc'_k M$, the chain $(A_{ki})_{i \ge 1}$ becomes
stationary.  Furthermore, by Lemma D(2), it consists of nonzero terms.  Therefore it
converges to a nonzero $S_{k-1}$-submodule $C_k
\subseteq \endosoc'_k M$.  From the definition of the
$A_{ki}$, we moreover derive 
$C_k = (J(S) C_{k+1}) \cap \endosoc'_k M$ for all $k \in \NN$, which yields a sequence of
inclusions 
$$C_1 \subseteq  J(S) C_2, \ C_2 \subseteq J(S) C_3,  \dots.$$

Since $C_1$ is an $S$-module, we can pick $i_1 \in I_1$ so that $M_{i_1} \cap C_1$ contains a
nonzero element
$x$.  The above string of inclusions then yields disjoint finite sets $L_n$ for
$n \in \NN$, together with families of homomorphisms $(f_{l_n})_{l_n \in L_n}$ in $J(S)$, each
of the form
$f_{l_n} = \pr_{v(l_n)} \cdot f_{l_n} \cdot
\inj_{u(l_n)} \in
\Hom_R(M_{u(l_n)}, M_{v(l_n)})$ for suitable indices $u(l_n) \in I_{n+1} \setminus I_n$,
$v(l_n) \in I_n$, and $v(l_1) = i_1$ for all $l_1 \in L_1$, such that
$$x \in \im \biggl( \sum_{l_1
\in L_1, \dots, l_n \in L_n} f_{l_1} f_{l_2}\cdots f_{l_n} \biggr)$$ for each $n$.  Here all of
the $f_{l_n}$'s are non-isomorphisms from $M_{u(l_n)}$ to $M_{v(l_n)}$, respectively, since
they belong to $J(S)$ (as usual, we identify homomorphisms $M_i \rightarrow M_j$ with elements
of $S$).  Note that the sets $\{u(l_n) \mid l_n \in L_n\}$ for $n \in \NN$ are pairwise
disjoint by construction; on the other hand, we permit repetitions
$u(l_k) = u(l_k')$ and $v(l_k) = v(l_k')$ for $l_k \ne l_k'$.  Next, we let $Y \subset
M_{i_1}$ be a maximal
$R$-submodule of
$M_{i_1}$ with $x \notin Y$, and consider the tree having as root the index $i_1$ and as
branches  all those concatenations
$f_{l_1} f_{l_2}\cdots f_{l_n}$ the image of which is not contained in $Y$.  By construction,
there is no upper bound on the lengths of the edge paths in this graph, and therefore K\"onig's
graph theorem guarantees the existence of a path of infinite length.  But, since $M_{i_1}/Y$
is a finitely cogenerated factor module of $M_{i_1}$, this means that the family
$(M_i)$ fails to be left semi-T-nilpotent.  Thus, our hypothesis ensures
finiteness of $I$.
\qed        
\enddemo

This smooths the road to

\proclaim{Theorem H}  Let $(M_i)_{i \in I}$ be a family of indecomposable
$R$-modules and $(M_t)_{t \in T}$ a trans\-versal of the isomorphism classes of the $M_i$. 
\smallskip

{\bf (I)}  If all of the $M_i$ have local endomorphism rings, the following statements are
equivalent:  

\quad {\rm (1)}  $\bigoplus_{i \in I} M_i$ is endo-noetherian, and
$(M_t)_{t \in T}$ is right semi-T-nilpotent.

\quad {\rm (2)}  $T$ is finite, and each $M_t$ is endo-noetherian. 
\smallskip

{\bf (II)}  The following statements are equivalent:

\quad {\rm (1)} $\bigoplus_{i \in I} M_i$ is endo-artinian, and $(M_t)_{t \in T}$ is left
semi-T-nilpotent.

\quad {\rm (2)}  $T$ is finite, and each $M_t$ is endo-artinian. 
\smallskip

{\bf (III)}  Suppose that $(M_t)_{t \in T}$ is a left semi-T-nilpotent and right T-nilpotent
family satisfying the following additional finiteness condition: $\ (\bullet)$ For any cofinite
subset $T' \subseteq T$, the endosocle of
$\bigoplus_{t \in T'} M_t$ is finitely generated.

Then $T$ is finite.
\endproclaim
 
\demo{Proof}  Start by observing that, in each of the first two parts, the implication
`(2)$\implies$(1)' is trivial.  Indeed, it suffices to note that both  endo-noetherianness and
endo-artinianness are inherited by arbitrary powers and finite direct sums.  Moreover, note
that, in proving the remaining implications, it is innocuous to assume that
$I=T$, since  the conditions involved are stable under passage from $\bigoplus_{i \in I} M_i$
to
$\bigoplus_{t \in T} M_t$.  The implications `(1)$\implies$(2)' under (I) and (II) now follow
from Lemma F and Lemma G, respectively;   note that (II.1) guarantees all of the $M_i$ to have
local endomorphsim rings, while this is a blanket hypothesis for (I).     

In order to prove (III), adopt all of the listed hypotheses.  By part (1) of Proposition E,
the $M_t$ again have local endomorphism rings.  Moreover, for each cofinite subset $T'
\subseteq T$, the direct sum
$\bigoplus_{t \in T'} M_t$ is finitely cogenerated over its endomorphism ring:  Indeed, in
view of Proposition E(2), right T-nilpotence of $(M_t)_{t \in T'}$ ensures that this direct
sum has essential endosocle and, by condition
$(\bullet)$, this endosocle is finitely generated.  Thus, once more, Lemma G yields finiteness
of $T$.     \qed         
\enddemo

The first part of the following corollary results from a combination of Theorem H(I,II) with
Proposition E; the second part is an immediate consequence of that proposition and Theorem
H(III).

\proclaim{Corollary I}  Let $(M_i)_{i \in I}$ be a family of indecomposable
$R$-modules and $(M_t)_{t \in T}$ a trans\-versal of the isomorphism classes of the $M_i$. 

(I)  Then the following statements are equivalent:

\quad {\rm (1)}   $\bigoplus_{i \in I} M_i$ is endo-noetherian, and
$(M_t)_{t \in T}$ is right T-nilpotent.

\quad {\rm (2)}  $T$ is finite, and each $M_t$ is endofinite. 
\smallskip

(II) Suppose that  $\bigoplus_{i \in I} M_i$ is $\Sigma$-algebraically compact, $(M_t)_{t \in
T}$ left semi-T-nilpotent, 

\quad\quad and 
\widestnumber\item{($\bullet$)\quad}\roster 
\item"($\bullet$)" for each cofinite subset $T' \subseteq T$, the endosocle of $\bigoplus_{t
\in T'} M_t$ is finitely generated.
\endroster

Then $T$ is finite. 
\qed\endproclaim

{\it Remarks.} None of the semi-T-nilpotence conditions in the various parts of Theorem H and
Corollary I are redundant, not even when all of the $M_i$ are finitely generated over a
finite dimensional algebra  --  see Examples O below.  

As for condition $(\bullet)$ in Theorem H(III) and Corollary I(II):  Dropping it means
jeopardizing the conclusion in either case, as is illustrated by Example J(a) below.  In case
the individual modules
$M_i$ have endosocles of finite length, condition $(\bullet)$ can clearly be weakened to the
requirement that, for each cofinite subset $T' \subseteq T$, the endosocle of $\bigoplus_{t
\in T'} M_t$ has finite support.  We do not know whether, in the above statements, condition
$(\bullet)$ can always be relaxed in this fashion, but point
out that finitely generated endosocles are not automatic in $\Sigma$-algebraically
compact modules (see Example J(b)).

\definition{Examples J} (a) The direct sum $M = \bigoplus_{p \ \text{prime}} \ZZ(p^{\infty})$
is a
$\Sigma$-algebraically compact $\ZZ$-module having as endomorphism ring the direct product of
the rings of $p$-adic integers, where $p$ runs through all primes.  Moreover, the family
$\bigl(\ZZ(p^{\infty})\bigr)_{p \ \text{prime}}$ is clearly left semi-T-nilpotent and right
T-nilpotent.  However, the endosocle of $M$ coincides with the $\ZZ$-socle, and thus has
infinite support.  This shows that the endosocle condition in Theorem H(III) and Corollary
I(II) is crucial.

(b)  Set $R = K[X_n \mid n \in \NN]/(X_n \mid n \in \NN)^k$, where $K$ is a field and $k \ge
2$.  Then the regular $R$-module is $\Sigma$-algebraically compact by \cite{\HuiZimI}. 
On the other hand, its endosocle,
$\endosoc R = \soc_R R = J^{k-1}$, fails to be finitely generated.  
\enddefinition

Prime targets for applications of the following consequence of Theorem H are families of
generic or finitely generated modules over Artin algebras.  In light of the subsequent
proposition addressing families of finitely generated modules over Artin algebras, the picture
can be sharpened in that situation. 

\proclaim{Corollary K}  For any family $(M_i)_{i \in I}$ of indecomposable endofinite
modules the following statements are equivalent:

{\rm (1)} The number of isomorphism types of $M_i$ is finite.

{\rm (2)}  $\bigoplus_{i \in I} M_i$ is endo-noetherian, and the family
$(M_i)$ is right semi-T-nilpotent.

{\rm (3)}  $\bigoplus_{i \in I} M_i$ is endo-artinian, and the family
$(M_i)$ is left semi-T-nilpotent.

{\rm (4)}  If $(M_t)_{t \in T}$ is a transversal of the isomorphism types of the $M_i$, then
$(M_t)$ is left and right semi-T-nilpotent and, for any cofinite subset $T' \subseteq T$, the
endosocle of $\bigoplus_{t \in T'} M_t$ has finite support in $T$.
\endproclaim

\demo{Proof}  In view of Theorem H, only the implication `(4)$\implies$(1)' requires some
justification:  Since the $M_i$ are endofinite, each singleton family
$\{M_i\}$ is right T-nilpotent by Proposition E(2), whence right T-nilpotence of the family
$(M_i)$ follows from right semi-T-nilpotence.     
\qed\enddemo

The next result exhibits a simplification of the picture for Artin
algebras.  Namely, a direct sum
$\bigoplus_{i \in I} M_i$ of finitely generated modules $M_i$ is $\Sigma$-algebraically compact
precisely when it is endo-artinian.  In other words, for Artin algebras, Proposition L
provides an `endo-artinian counterpart' to the final observation of Section 1: This
observation remains true on replacing `endo-noetherian' by `endo-artinian', and substituting
the ACC for (finite) matrix subgroups by the corresponding DCC.     

Throughout our excursion to the narrowed setting,
$\la$ will be an Artin algebra with center $C$, and $D: \lamod \rightarrow
\modla$ will denote the standard duality from the category of finitely generated left to the
category of finitely generated right $\la$-modules, i.e.,
$D = \Hom_C(-, E)$, where $E$ is the $C$-injective envelope of $C/J(C)$.  Note that
$D$ induces a Morita duality on $\Cmod$ as well; this obvious fact is useful,
because endo-submodules of $\la$-modules are of course
$C$-submodules.

\proclaim{Proposition L} If $\la$ is an Artin algebra and $(M_i)_{i \in I}$ a family of
finitely generated left $\la$-modules, the following conditions are equivalent:

{\rm (1)} $\bigoplus_{i \in I} M_i$ is $\Sigma$-algebraically compact.

{\rm (2)} $\bigoplus_{i \in I} M_i$ is endo-artinian.

{\rm (3)} $\bigoplus_{i \in I} D(M_i)$ is endo-noetherian.
\endproclaim

\demo{Proof} To verify `(1)$\implies$(3)', we adopt (1).  This means that $M = \bigoplus_{i \in
I} M_i$ satisfies the DCC for finite matrix subgroups, whence, by Proposition 3 of
\cite{\HuiZimII}, the dual $D(M) = \prod_{i \in I} D(M_i)$ has the ACC for finite
matrix subgroups.  In view of the final observation of Section 1, this entails (3).

To prove `(3)$\implies$(2)', we assume that $\bigoplus_{i \in I} D(M_i)$ is endo-noetherian. 
We begin by noting that a $C$-submodule $U$ of $\bigoplus_{i
\in I} M_i$ is an endo-submodule of $\bigoplus_{i \in I} M_i$ if and only if $U =
\bigoplus_{i \in I} U_i$, where the $U_i$ are subgroups of the
$M_i$ with the property that, for all $i,j \in I$ and all $f \in
\Hom_{\la}(M_i,M_j)$, the image $f(U_i)$ is contained in $U_j$. Now suppose that
$U$ is an endo-submodule of $\bigoplus_{i \in I} M_i$, and, for
$i \in I$, let
$\iota_i: U_i \hookrightarrow M_i$ be the canonical embedding.  Moreover, let $V_i
\subseteq D(M_i)$ denote the kernel of the dual map $D(\iota_i)$; i.e., $V_i = \{\rho
\in \Hom_C(M_i,E) \mid \rho|_{U_i} = 0\}$.  We claim that $V = \bigoplus_{i \in I} V_i$ is an
endo-submodule of
$\bigoplus_{i \in I} D(M_i)$. Clearly, $V$ is a $C$-submodule, whence it suffices to check
stability under homomorphisms $g \in \Hom_{\la}(D(M_i),D(M_j))$.  Suppose
$g = D(f)$ with $f \in \Hom_{\la}(M_j,M_i)$.  Dualizing the fact that $f\cdot
\iota_j = \iota_i \cdot (f|_{U_j})$, we obtain, for each $\rho \in V_i$, the equality
$g(\rho)|_{U_j} = (\rho\cdot f)|_{U_j} = (\rho|_{U_i}) \cdot (f|_{U_j}) = 0$; in other words,
$g(\rho)$ belongs to $V_j$ as required.

Let $U_1 \supseteq U_2 \supseteq U_3 \supseteq \cdots$ be a descending chain of
endo-submodules of $\bigoplus_{i \in I} M_i$, and label the embedding
$U_{k+1} \hookrightarrow U_k$ by $\phi_k$.  Extending the above notation, we moreover write
the embedding $U_{ki} = M_i \cap U_k \hookrightarrow M_i$ as
$\iota_{ki}$, and the kernel of the $C$-module epimorphism $D(\iota_{ki})$ as $V_{ki}$.  The
preceding paragraph shows that each of the direct sums $V_k = \bigoplus_{i \in I} V_{ki}$ is
in fact an endo-submodule of $\bigoplus_{i \in I} D(M_i)$.  Clearly, the
$V_k$ form an ascending chain which, by our hypothesis, becomes stationary.  This means that
the ascending chains $V_{1i} \subset V_{2i} \subset \cdots$ become uniformly
stationary on $I$, and, in view of the exact sequences
$$0 \rightarrow V_{ki} \rightarrow D(M_i) \rightarrow D(U_{ki}) \rightarrow 0,$$ we see that it
is again uniformly on $I$ that the maps $D(\phi_k |_{U_{ki}})$  --  and hence also the $\phi_k
|_{U_{ki}}$  --  become isomorphisms.  This proves our initial descending chain $(U_k)_{k \in
\NN}$ to become stationary.    

In view of the characterization of
$\Sigma$-algebraic compactness quoted in Section 1, the implication `(2)$\implies$(1)' is
clear.
$\qed$
\enddemo
  
The next corollary results from a
combination of points K and L with the Harada-Sai Lemma
\cite{\HarSai}.

\proclaim{Corollary M}  Suppose $\la$ is an Artin algebra and $(M_i)_{i
\in I}$  a family of left
$\la$-modules having uniformly bounded (finite) composition lengths.  Moreover, let $(M'_t)_{t
\in T}$ be a transversal of the isomorphism types of the indecomposable direct summands of the
$M_i$.  Then the following conditions are equivalent:

{\rm (1)} $T$ is finite.

{\rm (2)} $\bigoplus_{i \in I} M_i$ is endo-noetherian.

{\rm (3)} $\bigoplus_{i \in I} M_i$ is endo-artinian.

{\rm (4)} $\bigoplus_{i \in I} M_i$ is $\Sigma$-algebraically compact.

{\rm (5)}  For any cofinite subset $T' \subseteq T$, the endosocle of $\bigoplus_{t \in
T'}M'_t$ has finite support.  
\endproclaim

\demo{Proof}  It suffices to observe that, due to the Harada-Sai Lemma, our blanket hypothesis
forces the family $(M'_t)$ to be both left and right T-nilpotent. \qed \enddemo

The final corollary in this series highlights the role played by the endosocles of direct sums
of indecomposable representations.  It is an immediate consequence of Corollary M and the
(confirmed) second Brauer-Thrall conjecture.    

\proclaim{Corollary N}  For any finite dimensional algebra $\la$ over an algebraically closed
field, the following conditions are equivalent:

{\rm (1)}  $\la$ has finite representation type.

{\rm (2)}  For every family $(M_i)_{i \in I}$ of pairwise non-isomorphic indecomposable
modules of fixed finite dimension, the endosocle of $\bigoplus_{i \in I} M_i$ has finite
support. \qed  
\endproclaim

The following examples attest to the fact that none of the T-nilpotence conditions in Theorem H
and Corollaries I, K can be dropped without penalty.

\definition{Examples O} Let $\la = K\Gamma$ be the Kronecker algebra over a field $K$, where
$\Gamma$ is the quiver

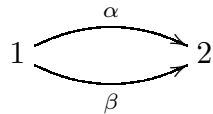
\ignore{
$$\xymatrixcolsep{4pc}
\xy\xymatrix{ 1 \ar@/^0.75pc/[r]^\alpha \ar@/_0.75pc/[r]_\beta &2 }\endxy$$ } 

(1) Moreover, let $(M_n)_{n \in \NN}$ be the family of all preinjective left $\la$-modules,
i.e., $M_1 = \la e_1/J(\la)e_1$, and for $n \ge 2$,
$M_n$ is the string module with graph

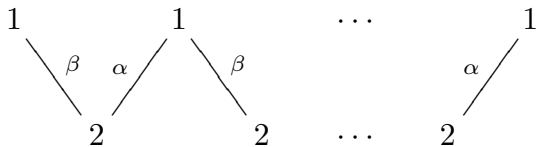
\ignore{
$$\xymatrixcolsep{1.25pc}
\xy\xymatrix{
 &1  \edge[dr]^{\beta} &&1 \edge[dl]_{\alpha}
\edge[dr]^{\beta} &&\cdots &&1 \edge[dl]_{\alpha}
 \\   &&2 &&2 &\cdots &2  }\endxy$$ }

\noindent having $n$ summands $S_1$ in its top.  Then the direct sum $\bigoplus_{n \in \NN}
M_n$ is $\Sigma$-algebraically compact by \cite{\Len, Theorem 4.6}, and Proposition L
guarantees that $\bigoplus_{n \in \NN} M_n$ is even endo-artinian.

(2)  Let $\la$ still be the Kronecker algebra, but this time consider the family $(N_n)_{n \in
\NN}$ of all preprojective left $\la$-modules, i\.e., $N_n = D(M_n')$, where $M_n'$ is the
preinjective right $\la$-module with top dimension $n$.  Since $\bigoplus_{n \in \NN} M_n'$ is
endo-artinian (see (1)), the direct sum $\bigoplus_{n \in \NN} N_n$ of the duals of the
$M_n'$ is endo-noetherian, again by Proposition L.   
\enddefinition

At this point, we return to an arbitrary base ring.  In light of our remarks concerning the
exchange property of direct sums (in Section 3), we finally obtain the following consequence
of Theorem H.  

\proclaim{Corollary P}  An endo-noetherian direct sum $\bigoplus_{i \in I} M_i$ of
indecomposable modules $M_i$ (over any ring $R$) has the exchange property if and only if each
$M_i$ has local endomorphism ring and the $M_i$ fall into finitely many isomorphism classes.
\qed
\endproclaim

Viewing the Kronecker examples in light of this fact, we see that the direct sum of the
preprojective modules over the Kronecker algebra fails to have the exchange property.  On the
other hand, the direct sum ranging over the preinjective modules does enjoy this property, as
do all algebraically compact modules (see
\cite{\HuiZimexchange}).

Variants of the ACC and DCC for endo-submodules have made numerous  appearances in the
representation theory of Artin algebras.  Our concluding remarks, relating ascending
endo-chain conditions of direct sums to preprojective partitions and left almost split
morphisms, are merely the `tip of the iceberg'.  (Naturally, certain relaxed descending
endo-chain conditions are linked up with
preinjective partitions and the existence of right almost split maps in a dual fashion.)
\smallskip

\noindent {\it Remarks.} 

Due to \cite{\AuSm} and \cite{\Huistrong}, respectively, we have the following connections
between weakened endo-noetherian conditions for direct sums $\bigoplus_{i \in I} M_i$ of
finitely generated indecomposable modules on one hand and the existence of preprojective
partitions (resp. strong preprojective partitions) on the other:  From now on, suppose
that
$(M_i)_{i \in I}$ is a family of finitely generated indecomposable modules, and $(M_t)_{t \in
T}$ a transversal of its isomorphism types.  Recall from \cite{\AuSm} that a {\it preprojective
partition} of
$(M_i)_{i \in I}$ is a partition 
$$\{M_t \mid t \in T \} = \bigsqcup_{n \le \omega} \Cal C_n$$ such that, for all $n < \omega$,
the set $\Cal C_n$ is a minimal finite generating set for $\bigsqcup_{n \le m \le \omega} \Cal
C_m$.  A {\it strong preprojective partition} extends this principle to arbitrary ordinals,
i.e., it is a
 partition  $\{M_t \mid t \in T \} = \bigsqcup_{\alpha \le \gamma} \Cal C_{\alpha}$, where
$\gamma$ is an ordinal number and each $\Cal C_{\alpha}$ is a minimal finite generating set
for $\bigsqcup_{\alpha \le \beta \le \gamma} \Cal C_{\beta}$.  If the $M_i$ have
perfect endomorphism rings, such a (strong) preprojective partition is unique in case of
existence.  The existence of a preprojective partition of
$(M_i)$, on the other hand, is equivalent to the following endo-condition for
$\bigoplus_{i \in I} M_i$:  For each cofinite subset $T_1 \subseteq T$, there exists a
cofinite subset $T_2 \subseteq T_1$ such that 
$$\bigoplus_{t \in T_1} M_t \subseteq \bigl( \End_R(\bigoplus_{t \in T} M_t) \bigr) \cdot 
\bigl(
\bigoplus_{t \in T_1 \setminus T_2} M_t \bigr).$$   Moreover, the condition that
$\bigoplus_{i \in I} M_i$ be endo-noetherian implies the existence of strong preprojective
partitions for arbitrary subfamilies of $(M_i)_{i \in I}$.  

The connection between endofiniteness conditions and strong preprojective partitions was
recently refined by Dung \cite{\Dung, 3.9 and 3.11} for the case where the $M_i$
have perfect endomorphism rings.  Namely, he established the following bridges among the
conditions labeled (a), (b), (c) below:

(a)  $\Hom_R(M_k, \bigoplus_{i \in I} M_i)$ is noetherian as a left module over
$\End_R(\bigoplus_{i \in I} M_i)$ for all $k\in I$. 

(b)  For every subfamily $\Cal C$ of $(M_i)$, the full subcategory $\add \Cal C$ of the
category of all finitely generated modules has left almost split morphisms. 

(c)  Every subfamily of $(M_i)$ has a strong preprojective partition.  
\smallskip

Under the given hypotheses for the $M_i$, condition (a) implies (b).  If, moreover, the
direct sum $\bigoplus_{i \in I} M_i$ is finitely generated over its endomorphism ring, then
(b) implies (c).

\enddocument